\documentclass[12pt]{amsart}
\usepackage{amsmath}
\usepackage{amsfonts,amssymb,amsthm}
\newtheorem{theorem}{Theorem}[section]
\newtheorem{lemma}[theorem]{Lemma}
\newtheorem{corollary}[theorem]{Corollary}
\theoremstyle{remark}
\newtheorem{remark}[theorem]{Remark}
\newtheorem{definition}[theorem]{Definition}
\def\vol{\mbox{vol}}
\def\tr{\mbox{tr}}
\def\diag{\mbox{diag}}
\begin{document}
\title{Some observations on the simplex}
\author{Igor Rivin}
\address{Mathematics Department, Temple University, 1805 N Broad St\\
  Philadelphia, PA 19122}
\address{Mathematics Department, Princeton University \\
Fine Hall,   Washington Rd \\
Princeton, NJ 08544} \email{rivin@math.temple.edu}
\thanks{The author was partially supported by the Department of
  Mathematical Sciences of the National Science Foundation. }
\keywords{convexity, simplex, volume,  symmetric, variation
  }
\begin{abstract}
We discuss some properties of the set of simplices in
$\mathbb{E}^n.$
\end{abstract}
\maketitle
\section{An inspirational example}
Consider triangles in the Euclidean plane $\mathbb{E}^2.$ The set
of congruence classes of triangles  is parametrized by their three
sidelengths (say $a,$ $b,$ $c$), while the set of similarity
classes of triangles is parametrized naturally by the triples of
angles $\alpha, \beta, \gamma.$

What can be said about the possible values of the triple $a, b,
c$? It is clearly necessary that the triangle inequalities be
satisfied:
\begin{gather*}
a  >  0, \  b  >  0,\  c  >  0,\\
a  <  b + c,\  b  <  a + c,\  c  <  a + b.
\end{gather*}
Since, given a triple $T=(a, b, c)$ satisfying the three triangle
inequalities we can easily construct the triangle with the
sidelengths prescribed by $T$, it follows that, as parametrized by
the lengths of sides, the set of non-degenerate triangles in
$\mathbb{E}^2$ is a \emph{convex cone}. A similar result holds for
the angles: a triple $\alpha, \beta, \gamma$ of \emph{positive}
numbers gives the angles of a non-degenerate triangle if and only
if
\begin{equation}
\alpha + \beta + \gamma = \pi.
\end{equation}
One can ask whether results of this type extend to simplices of
higher dimension, and we shall attempt to give as complete an
answer as possible.
\section{Some negative results}
First, we shall see that things appear to become less simple.
Firstly
\begin{theorem}
\label{nontri} Let $l_{0,1}, \dots, l_{n-1, n}$ be positive
numbers. In order for these to be the edge-lengths of a
non-degenerate simplex $v_0, \dots, v_n$ in $\mathbb{E}^n$, so
that $l_{ij} = d(v_i, v_j)$ it is necessary, but \emph{not}
sufficient that for any $v_i, v_j, v_k$ the lengths $l_{i,j},
l_{i, k}, l_{j, k}$ satisfy the requisite triangle inequalities.
\end{theorem}

\begin{proof}
The necessity is obvious. To prove the insufficiency, we construct
an example as follows: Let $n = 3,$ and let
$$l_{1,2} = l_{1, 3} = l_{2, 3} = 1,$$ while
$$l_{0, 1} = l_{0, 2} = l_{0, 3} = 1/2 + \epsilon.$$
The reader will verify that for any $\epsilon > 0$ all possible
triangle inequalities are satifsfied, whereas for a sufficiently
small $\epsilon$ no simplex with prescribed edgelengths exists.
\end{proof}

Things are even worse than that:

\begin{theorem}
Let $l_{i,j}, \quad 0\leq i, j \leq n$ be as in the statement of
Theorem \ref{nontri}. The set of simplices in $\mathbb{E}^n$ (for
$n > 2$) is \emph{not} convex when parametrized by the $l_{i,j}.$
\end{theorem}

\begin{proof} (The example is due to Peter Frankel of
  Budapest). Consider two simplices $A$ and $B$ as follows:
$$l_{0,1}(A) = l_{0, 2}(A) = l_{0, 3}(A) = l_{0, 1}(B) = l_{0, 2}(B) =
  l_{0, 3}(B) = 1,$$
while
\begin{eqnarray*}
l_{1, 2}(A) &=& \epsilon,\\
l_{2, 3}(A) &=& \sqrt{2},\\
l_{3, 1}(A) &=& \sqrt{2},\\
l_{1, 2}(B) &=& \sqrt{2},\\
l_{2, 3}(B) &=& \sqrt{2},\\
l_{3, 1}(B) &=& \epsilon.
\end{eqnarray*}
It is a simple geometric exercise to show that $A$ and $B$ really
and truly exist (for any $\epsilon > 0$) while, for $\epsilon$
small enough, there is no simplex $C$ with $l_{i,j}(C) =
l_{i,j}(A) + l_{i,
  j}(B), \quad \forall 0 \leq i<j \leq 3.$
(alternatively, the reader can peek ahead for a hint to a
non-geometric proof of non-existence of $C$.)
\end{proof}

We shall have to resort to less geometric methods to see what is
really going on.

\section{Linear algebra to the rescue}
Let $\Delta$ be a simplex with vertices $v_0 = \mathbf{0}, v_1,
\dots, v_n$ in $\mathbb{E}^n.$ The volume of the simplex $\Delta$
is non-zero if and only if the vectors $v_1, \dots, v_n$ are
linearly independent; in fact, its volume is given by
\begin{equation}
\label{volexp} \vol \Delta = \frac{1}{n!} \det V,
\end{equation}
where $V$ is the matrix whose columns are the vectors $v_1, \dots,
v_n.$

We now write the so-called \emph{Gram} matrix of $\Delta$:
$$G(\Delta) = V^t V.$$
Equivalently,
$$G_{ij}(\Delta) = \langle v_i, v_j \rangle.$$

\begin{theorem}
\label{gramchar} An $n \times n$ matrix $G$ is the Gram matrix of
a simplex $\Delta$ in $\mathbb{E}^n$ if and only if $G$ is
symmetric and positive definite.
\end{theorem}

\begin{proof}
First we prove the ``if'' direction. The symmetry of the inner
product in $\mathbb{E}^n$ implies that $G$ is symmetric
(alternately: $\left(V^t V\right)^t = V^t V.$) And for any vector
$x \in \mathbb{E}^n$:
$$x^t G(\Delta) x = x^t V^t V x = \langle V x, V x\rangle = \left\| V
x \right\| > 0,$$ if $V$ is non-singular and $x$ is non-zero,
hence positive definiteness.

To show the ``only if'' direction, let $G$ be a positive-definite
symmetric metrix. That implies that there is an orthogonal matrix
$S$ and a diagonal matrix $D$ all of whose diagonal elements are
positive, such that $G = S^t D S.$ Write $F = \sqrt{D}$ (the
elements of $F$ will be simply the non-negative  square roots of
the corresponding elements of $D$). Then, it is easy to see that
$G = S^t D S = S^t F^t F S = (F S)^t (FS).$ Setting $V = FS$ we
obtain our simplex.
\end{proof}

The fact, used in the proof of the above theorem, that a symmetric
matrix has an orthogonal matrix of eigenvectors is essentially
equivalent to the following
\begin{theorem}[Rayleigh-Ritz characterization]
\label{rr} The smallest
  eigenvalue $\lambda_0$ of a symmetric matrix $A$ can be written as
$$\lambda_0 = \min_{\|x\| = 1} \langle A x, x\rangle.$$
\end{theorem}
from which we have
\begin{corollary}
\label{evconc}
$$\lambda_0(A+B) \geq \lambda_0(A) + \lambda_0(B)$$ (in other words,
$\lambda_0$ is a concave function on the set of symmetric
matrices.
\end{corollary}
\begin{proof}
By the Rayleigh-Ritz characterization,
$$\lambda_0(A+B) = \min_{\|x\| = 1} \langle (A+B) x, x\rangle.$$
Denote the unit vector achieving the minimum by $x_{A+B}.$ Since
$$\langle A x_{A+B}, x_{A+B} \rangle \geq
\min_{\|x\| = 1} \langle A x, x\rangle,$$ and similarly
$$\langle B x_{A+B}, x_{A+B} \rangle \geq
\min_{\|x\| = 1} \langle B x, x\rangle,$$ the assertion follows.
\end{proof}
And the final corollary:
\begin{corollary}
\label{pdconv} The set of positive-definite symmetric matrices is
a \emph{convex
  cone} in the space of $n \times n$ symmetric matrices (which can be
  be naturally identified with $\mathbb{R}^{\binom{n+1}{2}}.$
\end{corollary}
\begin{proof}
This is an immediate consequence of Corollary \ref{evconc}, since
a matrix is positive-definite if and only if all of its
eigenvalues are positive.
\end{proof}

Corollary \ref{pdconv} can be restated as saying that the set of
Gram matrices of simplices is a convex cone. This does not see
very useful at the moment, since the entries of the Gram matrices
are some obscure scalar products, but this can be easily rectified
by observing the
\begin{quotation}
\textbf{Polarization identity}:
$$\langle v, w\rangle = \frac{1}{2} \left(\|v \|^2 + \|w\|^2 -
\|v-w\|^2\right).$$
\end{quotation}
When applied to the Gram matrix, this shows that
$$G_{ij}(\Delta) = \frac{1}{2}\left(l_{i,0}^2 + l_{j, 0}^2 -
l_{i,j}^2\right),$$ where $l_{i,j}$ are the edge lengths as
before, and $l_{i,i} = 0$ for any $i.$ We see that the Gram matrix
is a linear (matrix valued) function of the \emph{squares} of the
edge lengths of the simplex $\Delta,$ and hence we can put all of
the above together to see that:
\begin{theorem}
The set of non-degenerate simplices in $\mathbb{E}^n$ is a convex
cone when parametrized by the squares of the edgelengths.
\end{theorem}
We have now recovered our convexity of the set of simplices,
which, in view of the gloomy Section 2 can already be viewed as a
success. But, as they say, Wait! There is MORE!
\begin{remark}
The following is related to the results of this section, and
follows from a beautiful theorem of I.~J.~Schoenberg
(\cite{schoenberg}): the set of simplices parametrized by the
$\log$s of the edge lengths is \emph{starshaped} with respect to
the origin (corresponding to the regular simplex).
\end{remark}

\section{Volume}
In this section we shall use the following simple (as the reader
will see momentarily) result:

\begin{theorem}
\label{logdetthm} The function $\log \det A$ is a concave function
on the cone of positive definite symmetric matrices.
\end{theorem}

\begin{proof}
Concavity is equivalent to concavity on all lines, which is in
turn equivalent to the statement that
\begin{equation}
\label{secondderiv} \left.\frac{d^2 \log \det (A + t B)}{d
t^2}\right|_{t=0} < 0,
\end{equation}
for any positive definite $A$ and any symmetric $B.$ To show the
inequality (\ref{secondderiv}) we have to be able to differentiate
the determinant (or its logarithm directly). A couple of ways of
doing that will be shown below, but for now we will accept the
following as a fact of life:

\begin{equation}
\label{dld} \frac{d \log \det (A+t B)} {d t} = \tr ((A+tB)^{-1}
B).
\end{equation}
From which it follows easily (by linearity of trace) that
\begin{equation}
\label{ddld} \left.\frac{d^2 \log \det(A+t B)}{d t^2}\right|_{t=0}
= - \tr(A^{-1} B A^{-1} B).
\end{equation}
Now, since determinant is invariant under conjugation, we can
assume that $A$ is diagonal: $A = \diag(\lambda_1, \dots,
\lambda_n).$ Under that assumption, direct computation shows that
$$\tr(A^{-1}B A^{-1}B) = \sum_{0\leq i,j \leq n}
\frac{b_{ij}^2}{\lambda_i \lambda_j}.$$ Since all the summands are
nonnegative, and at least one is positive (since $B\neq 0$), the
theorem is proved.
\end{proof}
We finally obtain the following:
\begin{theorem}
\label{allvol} The volume of \emph{any fixed given face} (of any
dimension) of simplices parametrized by the squares of the
edgelengths is log-concave.
\end{theorem}
\begin{proof}
The result for the top-dimensional face follows from the
expression (\ref{volexp}) and Theorem \ref{logdetthm} above. Since
the Gram matrix of any given face is a submatrix of the Gram
matrix of the whole simplex, the result follows.
\end{proof}
Using the usual inequality theory (see \cite{gardner} for details)
we obtain the following result of a Brunn-Minkowsky type:
\begin{corollary}
\label{geommean} The $n$-th root of the volume of $n$-dimensional
simplices parametrized by the squares of the edgelengths is a
concave function.
\end{corollary}

These results can be used to produce hitherto unknown
characterizations of the regular simplex. First, an observation:

\begin{lemma}
\label{convlemma} Let $\Omega$ be a compact convex set in an
affine space, $G$ a group of automorphisms acting on $\Omega$, and
$f$ a concave function invariant under $G$ (that is, $f(G(x)) =
f(x), \quad \forall x \in \Omega$). Let $y$ be a point where $f$
achieves its maximum. Then $y$ is invariant under $G$ (that is,
$g(y) = y, \quad \forall g \in G.$)
\end{lemma}

\begin{proof}
Suppose that there is a $g \in G$ such that $z = g(y) \neq y.$
Since $\Omega$ is convex, the segment $S=[y, z]$ is contained in
$\Omega.$ The function $f$ is concave on $S,$ and so
$$
f\left(\frac{y+z}{2}\right) > \frac{f(y) + f(z)}{2} = f(y),
$$
contradicting the assumption that $y$ is a maximum of $f$.
\end{proof}

And now, the  applications:

\begin{theorem}
Let $\Delta_a$ be the set of $n$-dimensional simplices, such that
the sum of the squares of their edgelengths is equal to $a.$ Let
$$P_k = \prod_{\mbox{$k$-dimensional faces $f$ of $\Delta \in
\Delta_a$}} V(f).$$ Then $P_k$ is maximized when $\Delta$ is a
regular simplex (with edgelengths $2 a/(n+1)(n).$
\end{theorem}

\begin{theorem}
Let $\Delta_a$ be the set of $n$-dimensional simplices, such that
the sum of the squares of their edgelengths is equal to $a.$ Let
$$S_k = \sum_{\mbox{$k$-dimensional faces $f$ of $\Delta \in
\Delta_a$}} V^{1/k}(f).$$ Then $S_k$ is maximized when $\Delta$ is
a regular simplex (with edgelengths $2 a/(n+1)(n).$
\end{theorem}
\begin{proof}
The Theorems above follow immediately from Theorem \ref{allvol}
(Corollary \ref{geommean}, respectively) and Lemma
\ref{convlemma}.
\end{proof}

\section{How to differentiate the determinant}
In this section we demonstrate the truth of the formula
(\ref{dld}). Consider the family
$$M(t) = A + t B.$$
Clearly, $$\det M(t) = \det A \det(I + t A^{-1} B).$$ The
logarithmic derivative of $\det M(t)$ at $0$ is thus equal to the
derivative of $\det(I + t A^{-1} B)$ at 0. To simplify notation,
write $C = A^{-1} B.$ The $i$-th column of $I+t C$ equals
$\mathbf{e}_i + t \mathbf{c}_i,$ where $e_i$ is the $i$th standard
basis vector, and $\mathbf{c}_i$ is the $i$th column of $C.$ Using
the multilinearity of the determinant, we see that
\begin{eqnarray*}
\det(I + t C) &=& \\
 \det(\mathbf{e}_1 + t \mathbf{c}_1, \dots,
\mathbf{e}_n + t \mathbf{c}_n) &=& \\
1 + t \sum_{i=1}^n \det(\mathbf{e}_1, \dots, \mathbf{c}_i, \dots,
\mathbf{e}_n) + O(t^2) &=&\\ 1 + \tr(C) + O(t^2),
\end{eqnarray*}
showing that
$$\frac{d \det(1+t C)}{d t} = \tr C,$$
which completes the proof.
\section{Random facts about matrices}
\begin{definition}
Let $M$ be a matrix. The \emph{adjugate} $\widehat{M}$ of $M$ is
the matrix of cofactors of $M.$ That is, $\widehat{M}_{ij} =
(-1)^{i+j} \det M^{ij},$ where $M^{ij}$ is $M$ with the $i$-th row
and $j$-th column removed.
\end{definition}

The reason for this definition is
\begin{theorem}[Cramer's rule]
For any $n\times n$ matrix $M$ (over any commutative ring)
\begin{equation*}
M \widehat{M} = \widehat{M} M = (\det M) I(n),
\end{equation*}
where $I(n)$ is the $n \times n$ identity matrix.
\end{theorem}

We also need
\begin{definition}
The \emph{outer product} of vectors $v = (v_1, \dots, v_n)$ and $w
= (w_1, \dots, w_n)$ is the matrix $v \otimes w,$ defined as
follows:
\begin{equation*}
(v \otimes w)_{ij} = v_i w_j.
\end{equation*}
\end{definition}
\begin{remark}
As the notation suggests, the outer product is actually a tensor
prcduct, though it would be more correct to write $v \otimes w^*.$
The Dirac notation for the outer product would be $| v \rangle ||
\langle w|$ while the Dirac notation for the \emph{inner} product
would be $\langle v | w \rangle,$ this possibly explaining the
name \emph{outer product}.
\end{remark}
Consider an arbitrary vector $x = (x_1, \dots, x_n).$ We see that
\begin{equation}
\left[(v \otimes w) x\right]_k = \sum_{i=1}^n v_k w_i x_i = v_k
\langle w, x \rangle,
\end{equation}
so that
\begin{equation}
\label{outie} (v \otimes w) x = \langle w, x \rangle v.
\end{equation}
We see that $v \otimes w$ is a multiple of the projection operator
onto the subspace spanned by $v.$ In particular, in the case when
$\| v \| = 1,$ the operator $v \otimes v$ is the orthogonal
projection operator onto the subspace spanned by $v.$ Since $v
\otimes w$ is a rank $1$ operators all of its eigenvalues are
equal to $0.$ The one (potentially) nonzero eigenvalue equals
$\langle v, w\rangle.$

We now show:

\begin{theorem}
\label{nulladj} Suppose that $M$ has nullity $1$, and the null
space of $M$ is spanned by the vector $v,$ while the null space of
$M^t$ is spanned by the vector $w.$ Then
$$
\widehat{M}= c w \otimes v,
$$
\end{theorem}
\begin{proof}
Since $M$ is singular, we know that $\det M = 0,$ and so every
column of $\widehat M$ is in the null-space of $M.$ so, letting
$\mathbf{m}_i$ denote the $i$th column of $\widehat{M},$ we see
that
\begin{equation*}
\mathbf{m}_i = d_i v_i.
\end{equation*}
However, $\widehat{M^t} = \left(\widehat{M}\right)^t$  so
performing the computation on transposes we see that
\begin{equation*}
\mathbf{m^t}_i = e_i w_i.
\end{equation*}
We see that
\begin{equation*}
\widehat{M}_{ij} = d_i v_j = e_j w_i.
\end{equation*}
Writing $d_i = g_i w_i,$ and $e_j = h_j v_j,$ we see that, for
every pair $i, j,$ $g_i w_i v_j = h_j w_i v_j.$ The conclusion
follows.
\end{proof}

\begin{theorem}
The constant $c$ in the statement of the last theorem equals the
product of the nonzero eigenvalues of $M$ divided by the inner
product of $v$ and $w.$
\end{theorem}
\begin{proof}
By considering the characteristic polynomial of $M$ we see that
the product of the nonzero eigenvalues of $M$ equals the sum of
the principal $n-1$  minors. On the other hand, the principal
minors of $M$ equal the diagonal elements of $\widehat{M},$ so
\begin{equation*}
c \sum_{i=1}^n c w_i v_i = \prod_{i=1}^{n-1} \lambda_j.
\end{equation*}
\end{proof}
\begin{remark}
\label{regdet} By the discussion following Eq. \ref{outie}, the
product of nonzero eigenvalues of $M$ equals
$$\frac{\det \left(M + w \otimes v\right)}{\langle v, w\rangle}.$$
\end{remark}
\section{Back to simplices}
First we define the \emph{dual} Gram matrix of a simplex $S$.
First, let let $f_i$ be the unit outer normal to the $i$-th face.
Then
$$G_{ij}^* = \langle f_i, f_j \rangle.$$
In other words, the $ij$-th entry of the dual Gram matrix is the
cosine of the exterior dihedral angle between the $i$-th and the
$j$-th face.

\begin{lemma}
\label{cogram} The dual Gram matrix $G^*(\Delta)$ of a Euclidean
simplex $\Delta$ is symmetric and positive semi-definite, with
exactly one $0$ eigenvalue.
\end{lemma}
\begin{proof}
The proof proceeds exactly as the proof of Theorem \ref{gramchar}.
Any $n$ of the vectors $f_0, \dots, f_n$ are linearly independent,
and the corresponding $n \times n$ submatrix of the Gram matrix
$G^*$ is positive definite, whence the result.
\end{proof}
The above Lemma makes one wonder what the null space of $G^*$
might be. Happily, there is a complete answer, as follows:
\begin{theorem}
\label{nullthm} The null space of $G^*$ is generated by the vector
$\mathbf{A} = (A_0, \dots, A_n),$ where $A_i$ is the area of the
face $f_i.$
\end{theorem}
\begin{proof}
This follows from the divergence theorem, which states (in the
polyhedral case) that
\begin{equation}
\label{divergence} \sum_{i=0}^n A_i f_i =0.
\end{equation}
Taking the dot product of eq. (\ref{divergence}) with $f_j,$ we
see that
\begin{equation}
\sum_{i=0}^n A_i \langle f_i | f_j \rangle = 0,
\end{equation}
which is to say that the $j$-th coordinate of $G^* \mathbf{A}$
vanishes.
\end{proof}
\begin{corollary}
With notation as above,
\begin{equation*}
\frac{A_i^2}{A_j^2} = \frac{\widehat{M}_{ii}}{\widehat{M}_{jj}}
\end{equation*}
\end{corollary}
\begin{proof}
Immediate from Theorem \ref{nulladj} and Theorem \ref{nullthm}.
\end{proof}
\begin{remark}
This result has now been generalized to simplices in
$\mathbb{H}^3$ by the Novosibirsk State University Geometry
seminar.
\end{remark}
\bibliographystyle{alpha}

\end{document}